\pdfoutput=1
\RequirePackage{ifpdf}
\ifpdf % We are running pdfTeX in pdf mode
\documentclass[pdftex]{sigma}
\else
\documentclass{sigma}
\fi

\usepackage{easywncy}
\def\shp{\,\mathcyr{sh}\,}

\begin{document}

\renewcommand{\thefootnote}{$\star$}

\renewcommand{\PaperNumber}{061}

\FirstPageHeading

\ShortArticleName{The Algebra of a~$q$-Analogue of Multiple Harmonic Series}

\ArticleName{The Algebra of a~$\boldsymbol{q}$-Analogue\\
of Multiple Harmonic Series\footnote{This paper is a~contribution to the Special Issue in honor of Anatol
Kirillov and Tetsuji Miwa.
The full collection is available at \href{http://www.emis.de/journals/SIGMA/InfiniteAnalysis2013.html}
{http://www.emis.de/journals/SIGMA/InfiniteAnalysis2013.html}}}

\Author{Yoshihiro TAKEYAMA}

\AuthorNameForHeading{Y.~Takeyama}

\Address{Division of Mathematics, Faculty of Pure and Applied Sciences, University of Tsukuba,\\
Tsukuba, Ibaraki 305-8571, Japan}
\Email{\href{mailto:takeyama@math.tsukuba.ac.jp}{takeyama@math.tsukuba.ac.jp}}

\ArticleDates{Received June 27, 2013, in f\/inal form October 16, 2013; Published online October 22, 2013}

\Abstract{We introduce an algebra which describes the multiplication structure of a~family of $q$-series
containing a~$q$-analogue of multiple zeta values.
The double shuf\/f\/le relations are formulated in our framework.
They contain a~$q$-analogue of Hof\/fman's identity for multiple zeta values.
We also discuss the dimension of the space spanned by the linear relations realized in our algebra.}


\Classification{11M32; 33E20}

\vspace{-2mm}

\renewcommand{\thefootnote}{\arabic{footnote}}
\setcounter{footnote}{0}

\section{Introduction}

In this article we introduce an algebra to formalize the multiplication structure of a~$q$-analogue of
multiple zeta values.

An {\it admissible index} is an ordered set of positive integers $(k_{1}, \ldots, k_{r})$ with $k_{1} \ge
2$.
For an admissible index $\mathbf{k}=(k_{1}, \ldots, k_{r})$, the {\it multiple zeta value} (MZV)
$\zeta(\mathbf{k})$ is def\/ined by
\begin{gather*}
\zeta(\mathbf{k}):=\sum_{n_{1}>\cdots>n_{r}>0}\frac{1}{n_{1}^{k_{1}}\cdots n_{r}^{k_{r}}}.
\end{gather*}
The vector space spanned by MZVs over $\mathbb{Q}$ is closed under multiplication.
There are two ways to calculate the product of MZVs.
One way is to calculate the product directly from the above def\/inition of MZVs shuf\/f\/ling the indices
$n_{i}$.
Another way is to use an iterated integral representation, called the Drinfel'd integral~\cite{D, Zagier}.
By calculating the product of MZVs in two ways above, we obtain dif\/ferent expressions.
As a~result we get linear relations among MZVs, which are called the {\it double shuffle relations}.

In~\cite{H2} Hof\/fman gives an algebraic formulation to describe the multiplication structure of MZVs.
The two ways to calculate the product are realized as two dif\/ferent operations of multiplication on
a~non-commutative polynomial ring, which we call in this paper the harmonic product and the integral
shuf\/f\/le product, respectively.
Using the algebraic setup, an extension of the double shuf\/f\/le relations is given in~\cite{IKZ}, and it
is conjectured that it contains all linear relations among MZVs.

In this paper we consider the multiplication structure of a~$q$-analogue of MZVs.
Fix a~complex parameter $q$ such that $0<|q|<1$.
For an admissible index $\mathbf{k}=(k_{1}, \ldots, k_{r})$, a~{\it $q$-analogue of multiple zeta
values}~\cite{KKW, Z} is def\/ined by
\begin{gather}
\zeta_{q}(\mathbf{k}):=\sum_{n_{1}>\cdots>n_{r}>0}\frac{q^{(k_{1}-1)n_{1}+\cdots+(k_{r}-1)n_{r}}}{[n_{1}
]^{k_{1}}\cdots[n_{r}]^{k_{r}}},
\label{eq:def-qMZV}
\end{gather}
where
\begin{gather*}
[n]:=\frac{1-q^{n}}{1-q}
\end{gather*}
is the $q$-integer.
In the limit $q \to 1$, we restore the MZV $\zeta(\mathbf{k})$.

The harmonic product of $q$MZVs can be def\/ined naturally, and we also have an iterated integral
representation of $q$MZV~\cite{Z}.
However the vector space spanned by $q$MZVs over $\mathbb{Q}$ is not presumably closed under the
multiplication arising from the integral representation.
To overcome the dif\/f\/iculty we consider a~larger class of $q$-series allowing the factor $q^{n}/[n]$ in
the sum~\eqref{eq:def-qMZV}.
Such extension is proposed also in~\cite{OOZ}.
Then the enlarged vector space of $q$-series is closed under the harmonic product and the integral
shuf\/f\/le product.
The main result of this paper is to formulize the multiplication structure by extending Hof\/fman's algebra.
Thus we can consider the double shuf\/f\/le relations for $q$MZVs.

There are many linear relations over $\mathbb{Q}$ among $q$MZVs.
An important feature in the $q$-analogue case is that there are inhomogeneous linear relations in the
following sense.
For an admissible index $\mathbf{k}=(k_{1}, \ldots, k_{r})$, the {\it modified $q$MZV}
$\overline{\zeta}_{q}(\mathbf{k})$ is def\/ined by
\begin{gather}
\bar{\zeta}_{q}(\mathbf{k}):=(1-q)^{-\left|\mathbf{k}\right|}\zeta_{q}(\mathbf{k})=\sum_{n_{1}>\cdots>n_{r}
>0}\frac{q^{(k_{1}-1)n_{1}+\cdots+(k_{r}-1)n_{r}}}{(1-q^{n_{1}})^{k_{1}}\cdots(1-q^{n_{r}})^{k_{r}}},
\label{eq:def-qMZVbar}
\end{gather}
where $\left|\mathbf{k}\right|$ is the {\it weight} of $\mathbf{k}$ def\/ined by
$\left|\mathbf{k}\right|:=\sum\limits_{i=1}^{r}k_{i}$.
In~\cite{OT} it is observed that there are linear relations among the modif\/ied $q$MZVs with dif\/ferent
weight.
Taking the limit $q\to 1$ in such relations, the highest weight terms only survive and we obtain linear
relations for MZVs.
It suggests that we should consider the vector space spanned by the modif\/ied $q$MZVs rather than the
original $q$MZVs.

Our double shuf\/f\/le relations contain linear relations for the modif\/ied $q$MZVs.
However they do not suf\/f\/ice to get all linear relations.
In this article we also give some relations among $q$-series containing the factor $q^{n}/[n]$, which we
call the resummation duality, as a~supply of linear relations (see Theorem~\ref{thm:resummation-duality}
below).
By computer experiment it is checked that the double shuf\/f\/le relations and the resummation duality give
all linear relations among the modif\/ied $q$MZVs up to weight $7$.

The paper is organized as follows.
In Section~\ref{sec:double-shuffle} we give the algebraic setup to formalize the multiplication
structure of $q$MZVs.
To def\/ine the integral shuf\/f\/le product we make use of an extended version of a~$q$-analogue of
multiple polylogarithms (of one variable).
In Section~\ref{sec:dimension} we discuss the double shuf\/f\/le relations.
As an example we prove Hof\/fman's identity for $q$MZV (see Proposition~\ref{prop:Hoffman} below) in our
algebraic framework.
Note that it is derived from Ohno's relation and the duality for $q$MZV~\cite{B}.
At last we prove the resummation duality and show some computer experiment about the dimension of the
$\mathbb{Q}$-linear space spanned by the relations among the modif\/ied $q$MZVs obtained in this paper.

\section{Shuf\/f\/le products}
\label{sec:double-shuffle}

\subsection{Algebraic setup}

Let $\hbar$ be a~formal variable and $\mathcal{C}:=\mathbb{Q}[\hbar]$ the coef\/f\/icient ring.
Denote by $\mathfrak{H}$ the non-commutative polynomial algebra over $\mathcal{C}$ freely generated by
alphabet $\{x, y, \rho\}$.
Set
\begin{gather*}
\xi:=y-\rho,
\qquad
z_{k}:=x^{k-1}y,
\quad
k\ge1.
\end{gather*}
Let $\mathfrak{H}^{1}$ be the subalgebra of $\mathfrak{H}$ freely generated by the set $A:=\{\xi\}\cup
\{z_{k}\}_{k \ge 1}$.
Note that any element of $A$ is homogeneous and the degree of $\xi$ and $z_{k}$ $(k \ge 1)$ is $1$ and
$k$, respectively.

Def\/ine the $\mathcal{C}$-submodule $\widehat{\mathfrak{H}}^{0}$ of $\mathfrak{H}^{1}$ by
\begin{gather*}
\widehat{\mathfrak{H}}^{0}:=\mathcal{C}+\xi\mathfrak{H}^{1}+\sum_{k\ge2}z_{k}\mathfrak{H}^{1}.
\end{gather*}
We denote by $\mathfrak{H}^{0}$ the $\mathcal{C}$-submodule of $\widehat{\mathfrak{H}}^{0}$ generated by
$1$ and the words $z_{k_{1}} \dots z_{k_{r}}$ with $k_{1}\ge 2$ and $k_{2}, \ldots, k_{r} \ge 1$.

Hereafter we f\/ix a~complex parameter $q$ such that $0<|q|<1$.
We endow $\mathbb{C}$ with $\mathcal{C}$-module structure such that $\hbar$ acts as multiplication by $1-q$.
Denote by $\mathfrak{z}$ the $\mathcal{C}$-submodule of $\mathfrak{H}$ generated by $A$.
For a~positive integer $n$ we def\/ine the $\mathcal{C}$-linear map $I_{\cdot}(n): \mathfrak{z} \to
\mathbb{C}$ by
\begin{gather*}
I_{\xi}(n):=\frac{q^{n}}{[n]},
\qquad
I_{z_{k}}(n):=\frac{q^{(k-1)n}}{[n]^{k}}.
\end{gather*}
Note that
\begin{gather}
I_{\rho}(n)=I_{z_{1}-\xi}(n)=1-q.
\label{eq:I-rho}
\end{gather}

Now we def\/ine the $\mathcal{C}$-linear map $Z_{q}: \widehat{\mathfrak{H}}^{0} \to \mathbb{C}$ by $Z_{q}(1)=1$ and
\begin{gather*}
Z_{q}(u_{1}\dots u_{r}):=\sum_{n_{1}>\cdots>n_{r}>0}\;\prod_{i=1}^{r}I_{u_{i}}(n_{i}),
\end{gather*}
where $r \ge 1$ and $u_{i} \in A$.
The inf\/inite sum in the right hand side absolutely converges because there exists a~positive constant $M$
such that $|1/[n]|\le M$ for all $n\ge 1$.
If $\mathbf{k}=(k_{1}, \ldots, k_{r})$ is an admissible index, the value $Z_{q}(z_{k_{1}} \dots z_{k_{r}})$
is equal to the $q$MZV~\eqref{eq:def-qMZV}.

\subsection{Harmonic product}

We def\/ine the harmonic product on $\mathfrak{H}^{1}$ generalizing the algebraic formulation given
in~\cite{IKZ}.
Consider the commutative product $\circ$ on $\mathfrak{z}$ by setting
\begin{gather*}
z_{k}\circ z_{l}=z_{k+l}+\hbar z_{k+l-1},
\qquad
\xi\circ z_{k}=z_{k+1},
\qquad
\xi\circ\xi=z_{2}-\hbar\xi
\end{gather*}
for $k, l \ge 1$ and extending by $\mathcal{C}$-linearity.
Def\/ine the $\mathcal{C}$-bilinear product $*$ on $\mathfrak{H}^{1}$ inductively by setting
\begin{gather*}
1*w=w,
\qquad
w*1=w,
\\
(u_{1}w)*(u_{2}w')=u_{1}(w*u_{2}w')+u_{2}(u_{1}w*w')+(u_{1}\circ u_{2})(w*w')
\end{gather*}
for $w, w'\in\mathfrak{H}^{1}$ and $u_{1}, u_{2} \in A$.
It is commutative and associative because the product $\circ$ is commutative and associative.
Let us call $*$ the {\it harmonic product} on $\mathfrak{H}^{1}$.
Then the $\mathcal{C}$-submo\-du\-le~$\widehat{\mathfrak{H}}^{0}$ is a~subalgebra of $\mathfrak{H}^{1}$ with
respect to the harmonic product.
\begin{theorem}
\label{thm:harmonic-product}
For any $w, w'\in \widehat{\mathfrak{H}}^{0}$ we have $Z_{q}(w *w')=Z_{q}(w)Z_{q}(w')$.
\end{theorem}
\begin{proof}
For a~positive integer $N$ we def\/ine the $\mathcal{C}$-linear map $F_{\cdot}(N): \mathfrak{H}^{1} \to
\mathbb{C}$ by $F_{1}(N)=1$ and
\begin{gather}
F_{u_{1}\dots u_{r}}(N)=\sum_{N>n_{1}>\cdots>n_{r}>0}\;\prod_{i=1}^{r}I_{u_{i}}(n_{i}),
\label{eq:def-finite-harmonic-sum}
\end{gather}
where $u_{i} \in A$.
Note that $F_{u w}(N)=\sum\limits_{N>m>0}I_{u}(m)F_{w}(m)$ for $u \in A$ and $w \in \mathfrak{H}^{1}$.
We have $Z_{q}(w)=\lim\limits_{N \to \infty}F_{w}(N)$ for any $w \in \widehat{\mathfrak{H}}^{0}$, and hence it
suf\/f\/ices to prove that $F_{w *w'}(N)=F_{w}(N)F_{w'}(N)$ for words $w$, $w'$ in $A$ starting with $\xi$ or
$z_{k}$ $(k\ge 2)$.
Let us prove it by induction on the sum of the degrees of $w$ and $w'$.
Note that if $w=1$ or $w'=1$ the equality is trivial.
Let $w, w' \in \widehat{\mathfrak{H}}^{0}$ be words and $u_{1}, u_{2} \in A$.
Then
\begin{gather*}
F_{u_{1}w}(N)F_{u_{2}w'}(N)
=\sum_{N>m>0}I_{u_{1}}(m)F_{w}(m)F_{u_{2}w'}(m)+\sum_{N>m>0}I_{u_{2}}(m)F_{w'}(m)F_{u_{1}w}(m)
\\
\phantom{F_{u_{1}w}(N)F_{u_{2}w'}(N)=}
{}+\sum_{N>m>0}I_{u_{1}}(m)I_{u_{2}}(m)F_{w}(m)F_{w'}(m).
\end{gather*}
Now the desired equality follows from the induction hypothesis and
\begin{gather*}
I_{z_{k}}(m)I_{z_{l}}(m)=I_{z_{k+l}+\hbar z_{k+l-1}}(m),
\qquad
I_{\xi}(m)^{2}=I_{z_{2}-\hbar\xi}(m),
\qquad
I_{\xi}(m)I_{z_{k}}(m)=I_{z_{k+1}}(m)
\end{gather*}
for $k, l \ge 2$ and $m \ge 1$.
\end{proof}

\subsection{Integral shuf\/f\/le product}

Let us def\/ine the $\mathcal{C}$-bilinear product $\mathcyr{sh}$ on $\mathfrak{H}$
inductively as follows.
We set $1\shp w=w \shp 1=w$ for any $w \in \mathfrak{H}$.
For $u, v \in \{x, y, \rho\}$ and $w, w' \in \mathfrak{H}$, we set
\begin{gather*}
uw\shp vw'=u(w\shp vw')+v(uw\shp w')+\alpha(u,v)(w\shp w'),
\end{gather*}
where $\alpha(u, v)$ is determined by
\begin{gather*}
\alpha(x,x)=\hbar x,
\qquad
\alpha(x,y)=\alpha(y,x)=0,
\qquad
\alpha(y,y)=-y\rho
\end{gather*}
and
\begin{gather}
\alpha(u,\rho)=\alpha(\rho,u)=-u\rho
\label{eq:rho-prod}
\end{gather}
for $u \in \{x, y, \rho\}$.
Then the product $\mathcyr{sh}$ is commutative because of the symmetry of $\alpha$.
\begin{lemma}
\label{lem:rho-prod}
For any $w, w' \in \mathfrak{H}$, we have $\rho w\shp w'=\rho(w \shp w')$.
\end{lemma}
\begin{proof}
{}From the property~\eqref{eq:rho-prod}, we see that
\begin{gather*}
\rho w\shp u w'-\rho(w\shp u w')=u\left(\rho w\shp w'-\rho(w\shp w')\right)
\end{gather*}
for $u \in \{x, y, \rho\}$ and $w, w' \in \mathfrak{H}$.
Using this formula repeatedly we f\/ind that $\rho w\shp w'-\rho(w \shp
w')=w'\left(\rho w\shp 1-\rho(w \shp 1)\right)=0$ for $w, w' \in \mathfrak{H}$.
\end{proof}
\begin{proposition}
The product $\mathcyr{sh}$ is associative.
\end{proposition}
\begin{proof}
We prove $(w_{1}\shp w_{2})\shp w_{3}=w_{1}\shp (w_{2}\shp
w_{3})$ for $w_{i} \in \mathfrak{H}$ $(i=1, 2, 3)$ by induction on the sum of the degrees of $w_{1}$,
$w_{2}$ and $w_{3}$.
If $w_{i}=1$ for some $i$, it is trivial.
Suppose that $w_{i}=u_{i}w_{i}'$ $(i=1, 2, 3)$ for $u_{i}\in\{x, y, \rho\}$ and $w_{i}' \in \mathfrak{H}$.
Lemma~\ref{lem:rho-prod} implies that if $u_{i}=\rho$ for some $i$, the desired equality follows from the
induction hypothesis.
Thus we should check the associativity in the case where any $u_{i}$ is $x$ or $y$.
Here let us consider the case where $(u_{1}, u_{2}, u_{3})=(x, y, y)$.
We have
\begin{gather*}
(xw_{1}'\shp yw_{2}')\shp yw_{3}'=\big(x(w_{1}'\shp yw_{2}')+y(xw_{1}'\shp w_{2}')\big)\shp yw_{3}'
\\
\qquad{}
=x\big((w_{1}'\shp yw_{2}')\shp yw_{3}'\big)+y\big(x(w_{1}'\shp yw_{2}')\shp w_{3}'\big)
\\
\qquad\phantom{=}
{}+y\big((xw_{1}'\shp w_{2}')\shp yw_{3}'+y(xw_{1}'\shp w_{2}')\shp w_{3}'-\rho((xw_{1}'\shp w_{2}')\shp w_{3}')\big)
\\
\qquad
=x\big((w_{1}'\shp yw_{2}')\shp yw_{3}'\big)
\\
\qquad\phantom{=}
{}+y\left((xw_{1}'\shp w_{2}')\shp yw_{3}'+(xw_{1}'\shp yw_{2}')\shp w_{3}'-\rho((xw_{1}'\shp w_{2}')\shp w_{3}')\right).
\end{gather*}
Now apply the induction hypothesis and use the equality
\begin{gather*}
\rho\big(xw_{1}'\shp (w_{2}'\shp w_{3}')\big)=xw_{1}'\shp\big(\rho(w_{2}'\shp w_{3}')\big),
\end{gather*}
which follows from Lemma~\ref{lem:rho-prod}.
Then we obtain
\begin{gather*}
x\big(w_{1}'\shp (yw_{2}'\shp yw_{3}')\big)+y\big(xw_{1}'\shp \big(w_{2}'\shp yw_{3}'+yw_{2}'\shp w_{3}'-\rho(w_{2}'\shp w_{3}')\big)\big).
\end{gather*}
It is equal to $xw_{1}' \shp (yw_{2}'\shp yw_{3}')$.
The proof for the other cases is similar.
\end{proof}

Thus an associative commutative product $\mathcyr{sh}$ is def\/ined on $\mathfrak{H}$.
We call it the {\it integral shuffle product}.
\begin{proposition}
The $\mathcal{C}$-submodule $\widehat{\mathfrak{H}}^{0}$ of $\mathfrak{H}$ is closed under the integral
shuffle product.
\end{proposition}
\begin{proof}
First let us prove that
\begin{gather}
z_{k}w\shp z_{l}w'\in\sum_{j\ge\min{(k,l)}}z_{j}\mathfrak{H}^{1}
\label{eq:shuf-grading}
\end{gather}
for $k, l \ge 1$ and $w, w' \in \mathfrak{H}^{1}$.
If $k=l=1$ it follows from $\alpha(y, y)=-y\rho$.
If $k=1$ and $l\ge 2$, we f\/ind the above property by induction on $l$ using
\begin{gather*}
yw\shp z_{l}w'=y(w\shp z_{l}w')+x(yw\shp z_{l-1}w')
\end{gather*}
and $y=z_{1}$, $xz_{j}=z_{j+1}$ $(j \ge 1)$.
For $k, l \ge 2$ we obtain it by induction again from
\begin{gather*}
z_{k}w\shp z_{l}w'=x(z_{k-1}w\shp z_{l}w'+z_{k}w\shp z_{l-1}w'+\hbar z_{k-1}
w\shp z_{l-1}w').
\end{gather*}

It remains to prove that $\xi w\shp z_{k} w'$ and $\xi w \shp \xi w'$ belong to
$\widehat{\mathfrak{H}}^{0}$ for $w, w' \in \mathfrak{H}^{1}$ and $k \ge 2$.
It~follows from the property~\eqref{eq:shuf-grading} and
\begin{gather}
\xi w\shp z_{k}w'=\xi(w\shp z_{k}w')+x(yw\shp z_{k-1}w'),\label{eq:xishaxi}
\\
\xi w\shp \xi w'=\xi\big(w\shp \xi w'+\xi w\shp w'-\rho(w\shp w')\big).\tag*{\qed}
\end{gather}
\renewcommand{\qed}{}
\end{proof}

In the rest of this section we prove the following theorem.
\begin{theorem}
\label{thm:integral-shuffle-product}
For any $w, w'\in \widehat{\mathfrak{H}}^{0}$ we have $Z_{q}(w \shp w')=Z_{q}(w)Z_{q}(w')$.
\end{theorem}

Thus we def\/ine the two operations of multiplication, the harmonic product and the integral shuf\/f\/le
product, on $\widehat{\mathfrak{H}}^{0}$.
They describe the multiplication structure of a~family of $q$-series $Z_{q}(w)$ containing $q$MZVs.
Note that we can formally restore Hof\/fman's algebra for MZVs~\cite{H2} by setting $\hbar=0$ and $\rho=0$.

\subsection[A $q$-analogue of multiple polylogarithms]{A $\boldsymbol{q}$-analogue of multiple
polylogarithms}

To prove Theorem~\ref{thm:integral-shuffle-product} we introduce an extended version of a~$q$-analogue of
multiple polylogarithms (of one variable).
Denote by $\mathcal{F}$ the ring of holomorphic functions on the unit disk $|t|<1$.
We consider $\mathcal{F}$ as a~$\mathcal{C}$-module by $(\hbar f)(t):=(1-q)f(t)$ for $f \in \mathcal{F}$.
Def\/ine the $\mathcal{C}$-linear map $\widehat{\mathfrak{H}}^{0}\ni w \mapsto L_{w} \in \mathcal{F}$ by
$L_{1}(t)=1$ and
\begin{gather*}
L_{\xi w}(t):=\sum_{n=1}^{\infty}\frac{t^{n}}{[n]}F_{w}(n),
\qquad
L_{z_{k}w}(t):=\sum_{n=1}^{\infty}\frac{t^{n}}{[n]^{k}}F_{w}(n)
\end{gather*}
for $w \in \mathfrak{H}^{1}$ and $k\ge 2$, where $F_{w}(n)$ is def\/ined
by~\eqref{eq:def-finite-harmonic-sum}.

Consider the $q$-dif\/ference operator $\mathcal{D}_{q}$ def\/ined by
\begin{gather*}
(\mathcal{D}_{q}f)(t):=\frac{f(t)-f(qt)}{(1-q)t}.
\end{gather*}
To describe the function $\mathcal{D}_{q}L_{w}$ $(w \in \widehat{\mathfrak{H}}^{0})$ we introduce the two
maps $\Delta_{j}$ $(j=0, 1)$ as follows.
Set
\begin{gather*}
\widetilde{\mathfrak{H}}^{0}:=\xi\mathfrak{H}^{1}+\sum_{k\ge2}z_{k}\mathfrak{H}^{1},
\qquad
\mathfrak{h}^{\ge a}:=\mathcal{C}+\sum_{k\ge a}z_{k}\mathfrak{H}^{1},
\quad
a=1,2.
\end{gather*}
Let $\Delta_{0}: \mathfrak{h}^{\ge 2} \to \widetilde{\mathfrak{H}}^{0}$ be the $\mathcal{C}$-linear map
def\/ined by
\begin{gather*}
\Delta_{0}(1)=0,
\qquad
\Delta_{0}(z_{k}w)=
\begin{cases}
\xi w,&k=2,
\\
z_{k-1}w,&k\ge3
\end{cases}
\end{gather*}
for $w \in \mathfrak{H}^{1}$.
Next we def\/ine the $\mathcal{C}$-linear map $\Delta_{1}: \mathfrak{h}^{\ge 1} \to
\widehat{\mathfrak{H}}^{0}$ by $\Delta_{1}(1)=1$ and
\begin{gather*}
\Delta_{1}(z_{k}w)=\left(\sum_{a=2}^{k}\binom{k-1}{a-1}(-\hbar)^{k-a}z_{a}+(-\hbar)^{k-1}\xi\right)w
\end{gather*}
for $k \ge 1$ and $w \in \mathfrak{H}^{1}$.

Now note that $\widehat{\mathfrak{H}}^{0}$ is decomposed into the $\mathcal{C}$-submodules
\begin{gather*}
\widehat{\mathfrak{H}}^{0}=\mathfrak{h}^{\ge2}\oplus\left(\bigoplus_{r\ge0}\xi\rho^{r}\mathfrak{h}^{\ge1}
\right).
\end{gather*}
\begin{proposition}
\label{prop:multilog-q-difference}
For $w \in \mathfrak{h}^{\ge 2}$ we have
\begin{gather}
(\mathcal{D}_{q}L_{w})(t)=\frac{1}{t}L_{\Delta_{0}(w)}(t).
\label{eq:multilog-difference-0}
\end{gather}
For $w \in \mathfrak{h}^{\ge 1}$ and $r \ge 0$ it holds that
\begin{gather}
(\mathcal{D}_{q}L_{\xi\rho^{r}w})(t)=\frac{((1-q)t)^{r}}{(1-t)^{r+1}}L_{\Delta_{1}(w)}(t).
\label{eq:multilog-difference-1}
\end{gather}
\end{proposition}
\begin{proof}
The equality~\eqref{eq:multilog-difference-0} follows from $\mathcal{D}_{q}(t^{n})=[n]t^{n-1}$ for $n \ge 0$.
Let us prove~\eqref{eq:multilog-difference-1}.
If $w=u_{1}\dots u_{s} \in \mathfrak{h}^{\ge 1}$ $(u_{i} \in A)$ is a~word we have
\begin{gather*}
L_{\xi\rho^{r}w}(t)=(1-q)^{r}\sum_{n>n_{1}>\cdots>n_{s}>0}\binom{n-n_{1}-r}{r}\frac{t^{n}}{[n]}\prod_{i=1}^{s}I_{u_{i}}(n_{i})
\end{gather*}
because of~\eqref{eq:I-rho}.
Therefore
\begin{gather*}
\big(\mathcal{D}_{q}L_{\xi\rho^{r}w}\big)(t)=(1-q)^{r}\sum_{n_{1}>\cdots>n_{s}>0}\left(\sum_{n=n_{1}+1}
^{\infty}\binom{n-n_{1}-r}{r}t^{n-1}\right)\prod_{i=1}^{s}I_{u_{i}}(n_{i})
\\
\phantom{\big(\mathcal{D}_{q}L_{\xi\rho^{r}w}\big)(t)}
=\frac{((1-q)t)^{r}}{(1-t)^{r+1}}\sum_{n_{1}>\cdots>n_{s}>0}t^{n_{1}}\prod_{i=1}^{s}I_{u_{i}}(n_{i}).
\end{gather*}
Here we used the equality
\begin{gather}
\frac{1}{(1-x)^{k+1}}=\sum_{j=0}^{\infty}\binom{k+j}{j}x^{j},
\qquad
|x|<1,
\label{eq:expansion}
\end{gather}
for any non-negative integer $k$.
Now the equality~\eqref{eq:multilog-difference-1} follows from
\begin{gather*}
t^{n}I_{z_{k}}(n)=\sum_{a=1}^{k}\binom{k-1}{a-1}\big({-}(1-q)\big)^{k-a}\frac{t^{n}}{[n]^{a}}
\end{gather*}
for $k \ge 1$.
\end{proof}

\subsection[Multiplication structure of the $q$-analogue of multiple polylogarithms]{Multiplication
structure of the $\boldsymbol{q}$-analogue of multiple polylogarithms}

Let us prove that the set of functions $\{L_{w}\}_{w \in \widehat{\mathfrak{H}}^{0}}$ is closed under
multiplication.
Def\/ine the $\mathcal{C}$-linear map $I_{0}: \widetilde{\mathfrak{H}}^{0} \to \mathfrak{h}^{\ge
2}\cap\widetilde{\mathfrak{H}}^{0}$ by
\begin{gather*}
I_{0}(\xi w)=z_{2}w,
\qquad
I_{0}(z_{k}w)=z_{k+1}w,
\qquad
w\in\mathfrak{H}^{1}, \quad k\ge2,
\end{gather*}
and the $\mathcal{C}$-linear map $I_{1}: \widehat{\mathfrak{H}}^{0} \to \mathfrak{h}^{\ge 1}$ by
\begin{gather*}
I_{1}(1)=1,
\qquad
I_{1}(\xi w)=z_{1}w,
\qquad
I_{1}(z_{k}w)=\sum_{a=1}^{k}\binom{k-1}{a-1}\hbar^{k-a}z_{a}
\end{gather*}
for $w \in \mathfrak{H}^{1}$ and $k \ge 2$.
We have the following property.
\begin{lemma}
\label{lem:delta-varphi-identity}
The maps $\Delta_{0}I_{0}$ and $\Delta_{1}I_{1}$ are identities on $\widetilde{\mathfrak{H}}^{0}$ and
$\widehat{\mathfrak{H}}^{0}$, respectively.
\end{lemma}
\begin{proposition}\label{prop:q-integration}\qquad
\begin{enumerate}\itemsep=0pt
\item[$(1)$] Let $w \in \widetilde{\mathfrak{H}}^{0}$.
Suppose that $f \in \mathcal{F}$ satisfies $f(0)=0$ and $(\mathcal{D}_{q}f)(t)=L_{w}(t)/t$.
Then $f=L_{I_{0}(w)}$.

\item[$(2)$] Let $w \in \widehat{\mathfrak{H}}^{0}$ and $r \ge 0$.
Suppose that $f \in \mathcal{F}$ satisfies $f(0)=0$ and
\begin{gather*}
(\mathcal{D}_{q}f)(t)=\frac{((1-q)t)^{r}}{(1-t)^{r+1}}L_{w}(t).
\end{gather*}
Then $f=L_{\xi\rho^{r}I_{1}(w)}$.
\end{enumerate}
\end{proposition}
\begin{proof}
Note that the initial value problem $\mathcal{D}_{q}f=g$ and $f(0)=0$ for a~given $g \in \mathcal{F}$ has
a~unique solution in $\mathcal{F}$ if it exists.
Therefore it suf\/f\/ices to check that the function $f$ given above is a~solution to the initial value
problems in (1) or (2).
We have $f(0)=0$ because the image of $I_{0}$ or $\xi\rho^{r}I_{1}$ is contained in
$\widetilde{\mathfrak{H}}^{0}$.
Proposition~\ref{prop:multilog-q-difference} and Lemma~\ref{lem:delta-varphi-identity} imply that the
function $f$ is a~solution.
\end{proof}

To write down the structure of multiplication of the functions $L_{w}$ $(w \in
\widehat{\mathfrak{H}}^{0})$, let us def\/ine the commutative $\mathcal{C}$-bilinear product $\star$ on
$\widehat{\mathfrak{H}}^{0}$.
Set $1\star w=w \star 1=w$ for $w \in \widehat{\mathfrak{H}}^{0}$.
In general we def\/ine the product $\star$ inductively as follows.
For $w, w' \in \mathfrak{h}^{\ge 2}$ we set
\begin{gather*}
w\star w'=I_{0}\left(\Delta_{0}(w)\star w'+w\star\Delta_{0}(w')-\hbar\Delta_{0}(w)\star\Delta_{0}
(w')\right).
\end{gather*}
For $w \in \mathfrak{h}^{\ge 2}$, $w'\in\mathfrak{h}^{\ge 1}$ and $r \ge 0$, set
\begin{gather*}
w\star\xi\rho^{r}w'=I_{0}(\Delta_{0}(w)\star\xi\rho^{r}w')+\xi\rho^{r}I_{1}\big(\big(w-\hbar\Delta_{0}(w)\big)\star\Delta_{1}(w')\big).
\end{gather*}
For $w, w' \in \mathfrak{h}^{\ge 1}$ and $r, s\ge 0$,
\begin{gather*}
\xi\rho^{r}w\star\xi\rho^{s}w'=\xi\rho^{r}I_{1}(\Delta_{1}(w)\star\xi\rho^{s}w')+\xi\rho^{s}I_{1}
\big(\xi\rho^{r}w\star\Delta_{1}(w')\big)
\\
\phantom{\xi\rho^{r}w\star\xi\rho^{s}w'=}{}
-\xi\rho^{r+s+1}I_{1}\big(\Delta_{1}(w)\star\Delta_{1}(w')\big).
\end{gather*}
Since the image of $I_{0}$ is contained in $\widetilde{\mathfrak{H}}^{0}$, the product $\star$ is
well-def\/ined.
\begin{proposition}
\label{prop:multilog-product}
For any $w, w' \in \widehat{\mathfrak{H}}^{0}$ we have $L_{w \star w'}=L_{w}L_{w'}$.
\end{proposition}
\begin{proof}
It suf\/f\/ices to prove in the case where $w$ and $w'$ are homogeneous.
Let us prove it by induction on the sum of the degrees of~$w$ and~$w'$.
Note that the desired equality is trivial if~$w$ or~$w'$ belongs to $\mathcal{C}$.
Otherwise the function $L_{w}L_{w'}$ has a~zero at $t=0$.
Now calculate $\mathcal{D}_{q}(L_{w}L_{w'})$ by using the formula
\begin{gather*}
\big(\mathcal{D}_{q}(fg)\big)(t)=(\mathcal{D}_{q}f)(t)g(t)+f(t)(\mathcal{D}_{q}g)(t)-(1-q)t(\mathcal{D}_{q}
f)(t)(\mathcal{D}_{q}g)(t)
\end{gather*}
for $f, g \in \mathcal{F}$.
The $q$-dif\/ference of $L_{w}$ and $L_{w'}$ is written in terms of the maps $\Delta_{0}$ and $\Delta_{1}$
as described in Proposition~\ref{prop:multilog-q-difference}.
Here note that if $w$ is homogeneous the degree of $\Delta_{0}(w)$ is less than that of $w$.
Now the induction hypothesis implies that $\mathcal{D}_{q}(L_{w}L_{w'})$ is given in terms of the product~$\star$.
Use Proposition~\ref{prop:q-integration} to restore the original function $L_{w}L_{w'}$, and we get the
desired equality from the def\/inition of the product~$\star$.
\end{proof}

\subsection{Proof of Theorem~\ref{thm:integral-shuffle-product}}

Let us prove Theorem~\ref{thm:integral-shuffle-product}.
First we describe a~relation between the $q$MZV and the function $L_{w}$.
\begin{lemma}
\label{lem:evaluation}
Define the $\mathcal{C}$-linear map $e: \widehat{\mathfrak{H}}^{0} \to \widehat{\mathfrak{H}}^{0}$ by
setting $e(1)=1$ and
\begin{gather*}
e(\xi w)=\xi w,
\qquad
e(z_{k}w)=\left(\sum_{a=2}^{k}\binom{k-2}{a-2}\hbar^{k-a}z_{a}\right)w
\end{gather*}
for $w \in \mathfrak{H}^{1}$ and $k \ge 2$.
Then we have $L_{w}(q)=Z_{q}(e(w))$ for any $w \in \widehat{\mathfrak{H}}^{0}$.
\end{lemma}
\begin{proof}
It follows from $q^{n}/[n]=I_{\xi}(n)$ and $q^{n}/[n]^{k}=I_{e(z_{k})}(n)$ for $k \ge 2$ and $n \ge 1$.
\end{proof}

Note that the map $e$ given in Lemma~\ref{lem:evaluation} is an isomorphism on the $\mathcal{C}$-module
$\widehat{\mathfrak{H}}^{0}$.
Its inverse is given by $e^{-1}(1)=1$, $e^{-1}(\xi w)=\xi w$ and
\begin{gather*}%\label{eq:e-inverse}
e^{-1}(z_{k}w)=\left(\sum_{a=2}^{k}\binom{k-2}{a-2}(-\hbar)^{k-a}z_{a}\right)w
\end{gather*}
for $w \in \mathfrak{H}^{1}$ and $k \ge 2$.

Now Theorem~\ref{thm:integral-shuffle-product} is reduced to the following proposition because of
Proposition~\ref{prop:multilog-product} and Lemma~\ref{lem:evaluation}.
\begin{proposition}
\label{prop:star-shuffle}
It holds that $e\star(e^{-1}\times e^{-1})=\mathcyr{sh}$ on
$\widehat{\mathfrak{H}}^{0}\times\widehat{\mathfrak{H}}^{0}$.
\end{proposition}

In the proof of Proposition~\ref{prop:star-shuffle} we use the properties below.
\begin{lemma}\label{lem:formulas}   \qquad
\begin{enumerate}\itemsep=0pt
\item[$(1)$] $e(z_{k}w)=(x+\hbar)e(z_{k-1}w)$ for $k \ge 3$ and $w \in \mathfrak{H}^{1}$.

\item[$(2)$] $(1-\hbar \Delta_{0})e^{-1}=\Delta_{1}$ on $\mathfrak{h}^{\ge 2}$.

\item[$(3)$] $\Delta_{1}(z_{k}w)=-\hbar\Delta_{1}(z_{k-1}w)+e^{-1}(z_{k}w)$ for $k \ge 2$ and $w \in
\mathfrak{H}^{1}$.

\item[$(4)$] $\Delta_{0}e^{-1}(z_{k}w)=\Delta_{1}(z_{k-1}w)$ for $k \ge 2$ and $w \in \mathfrak{H}^{1}$.

\item[$(5)$] $eI_{0}e^{-1}(\xi w)=xyw$ for $w \in \mathfrak{H}^{1}$, and $eI_{0}e^{-1}(w)=(x+\hbar)w$ for $w \in
\sum\limits_{a \ge 2}z_{a}\mathfrak{H}^{1}$.

\item[$(6)$] $I_{1}I_{0}(w)=(x+\hbar)I_{1}(w)$ for $w \in \widehat{\mathfrak{H}}^{0}$.
\end{enumerate}
\end{lemma}

The proof is straightforward.
\begin{lemma}\label{lem:main}
For any $w, w' \in \mathfrak{h}^{\ge 1}$ it holds that
$I_{1}(\Delta_{1}(w)\star\Delta_{1}(w'))=w\shp w'$.
\end{lemma}
\begin{proof}
We can assume without loss of generality that $w$ and $w'$ are homogeneous.
If $w=1$ or $w'=1$, it is trivial since $I_{1}\Delta_{1}$ is the identity on $\mathfrak{h}^{\ge 1}$
(Lemma~\ref{lem:delta-varphi-identity}).
Let us prove the desired equality by induction on the sum of the degrees of $w$ and $w'$.

First consider the case where $w=z_{1}\rho^{r}w_{1}$ and $w'=z_{1}\rho^{s}w_{2}$ for $r, s\ge 0$ and
$w_{1}, w_{2} \in \mathfrak{H}^{1}$.
From the def\/inition of $\Delta_{1}$ and $\star$ we have
\begin{gather*}
I_{1}\big(\Delta_{1}(w)\star\Delta_{1}(w')\big)=I\big(\xi\rho^{r}w_{1}\star\xi\rho^{s}w_{2}\big)
\\
\hphantom{I_{1}\big(\Delta_{1}(w)\star\Delta_{1}(w')\big)}{}
=I_{1}\big(\xi\rho^{r}I_{1}(\Delta_{1}(w_{1})\star\xi\rho^{s}w_{2})
+\xi\rho^{s}I_{1}\big(\xi\rho^{s}w_{1}\star\Delta_{1}(w_{2})\big)\\
\hphantom{I_{1}\big(\Delta_{1}(w)\star\Delta_{1}(w')\big)=}{}
 -\xi\rho^{r+s+1}I_{1}\big(\Delta_{1}(w_{1})\star\Delta_{1}(w_{2})\big)\big)
\\
\hphantom{I_{1}\big(\Delta_{1}(w)\star\Delta_{1}(w')\big)}{}
=y\rho^{r}I_{1}\big(\Delta_{1}(w_{1})\star\Delta_{1}(y\rho^{s}w_{2})\big)
+y\rho^{s}I_{1}\big(\Delta_{1}(y\rho^{s}w_{1})\star\Delta_{1}(w_{2})\big)
\\
\hphantom{I_{1}\big(\Delta_{1}(w)\star\Delta_{1}(w')\big)=}{}
 -y\rho^{r+s+1}I_{1}\big(\Delta_{1}(w_{1})\star\Delta_{1}(w_{2})\big).
\end{gather*}
Apply the induction hypothesis and we get
\begin{gather*}
y\rho^{r}(w_{1}\shp y\rho^{s}w_{2})+y\rho^{s}(y\rho^{s}w_{1}\shp w_{2})-y\rho^{r+s+1}
(w_{1}\shp w_{2}).
\end{gather*}
It is equal to $z_{1}\rho^{r}w_{1}\shp z_{1}\rho^{s}w_{2}$ because of Lemma~\ref{lem:rho-prod}.

Next let us consider the case where $w=z_{1}\rho^{r}w_{1}$ and $w'=z_{k}w_{2}$ for $r \ge 0$, $k \ge 2$ and
$w_{1}, w_{2} \in \mathfrak{H}^{1}$.
Using Lemma~\ref{lem:formulas}~(3) and the induction hypothesis, we f\/ind that
\begin{gather*}
I_{1}\big(\Delta_{1}(w)\star\Delta_{1}(w')\big)=-\hbar\big(y\rho^{r}w_1\shp z_{k-1}w_2\big)
+I_{1}\big(\xi\rho^{r}w_1\star e^{-1}(z_{k}w_{2})\big).
\end{gather*}
From Lemma~\ref{lem:formulas}~(2),~(4) and $e^{-1}(z_{k}w) \in \mathfrak{h}^{\ge 2}$, we have
\begin{gather*}
\xi\rho^{r}w_1\star e^{-1}(z_{k}w_{2})
=I_{0}\big(\Delta_{1}(y\rho^{r}w_{1})\star\Delta_{1}(z_{k-1}w_{2})\big)+\xi\rho^{r}I_{1}\big(\Delta_{1}(w_{1})\star\Delta_{1}(z_{k}w_{2})\big).
\end{gather*}
Use Lemma~\ref{lem:formulas}~(6) to calculate the image of the f\/irst term by~$I_{1}$.
Now we can apply the induction hypothesis and see that
\begin{gather*}
I_{1}\big(\xi\rho^{r}w_1\star e^{-1}(z_{k}w_{2})\big)=(x+\hbar)(y\rho^{r}w_{1}\shp z_{k-1}w_{2}
)+y\rho^{r}(w_{1}\shp z_{k}w_{2}).
\end{gather*}
Therefore{\samepage
\begin{gather*}
I_{1}\big(\Delta_{1}(w)\star\Delta_{1}(w')\big)=x\big(y\rho^{r}w_{1}\shp z_{k-1}w_{2}\big)+y\rho^{r}(w_{1}
\shp z_{k}w_{2}).
\end{gather*}
It is equal to $z_{1}\rho^{r}w_{1}\shp z_{k}w_{2}$.}

Finally suppose that $w=z_{k}w_{1}$ and $w'=z_{l}w_{2}$ for $k, l \ge 2$ and $w_{1}, w_{2} \in
\mathfrak{H}^{1}$.
From Lemma~\ref{lem:formulas}~(3) we get
\begin{gather*}
\Delta_{1}(w)\star\Delta_{1}(w')
=e^{-1}(z_{k}w_{1})\star e^{-1}(z_{l}w_{2})-\hbar\Delta_{1}(z_{k}w_{1})\star\Delta_{1}(z_{l-1}w_{2})
\\
\phantom{\Delta_{1}(w)\star\Delta_{1}(w')=}
{}-\hbar\Delta_{1}(z_{k-1}w_{1})\star\Delta_{1}(z_{l}w_{2})
-\hbar^{2}\Delta_{1}(z_{k-1}w_{1})\star\Delta_{1}(z_{l-1}w_{2}).
\end{gather*}
Using the induction hypothesis we have
\begin{gather*}
I_{1}\big(\Delta_{1}(w)\star\Delta_{1}(w')\big)
=I_{1}\big(e^{-1}(z_{k}w_{1})\star e^{-1}(z_{l}w_{2})\big)-\hbar z_{k}w_{1}\shp z_{l-1}w_{2}
\\
\phantom{I_{1}\big(\Delta_{1}(w)\star\Delta_{1}(w')\big)=}
{}-\hbar z_{k-1}w_{1}\shp z_{l}w_{2}-\hbar^{2}z_{k-1}w_{1}\shp z_{l-1}w_{2}.
\end{gather*}
Since $e^{-1}(z_{k}w_{1})$ and $e^{-1}(z_{l}w_{2})$ belong to $\mathfrak{h}^{\ge 2}$, we see that
$I_{1}(e^{-1}(z_{k}w_{1})\star e^{-1}(z_{l}w_{2}))$ is equal to
\begin{gather*}
(x+\hbar)I_{1}\big(\Delta_{1}(z_{k-1}w_{1})\star\Delta_{1}(z_{l}w_{2})\,{+}\,\Delta_{1}(z_{k}w_{1}
)\star\Delta_{1}(z_{l-1}w_{2})\,{+}\,\hbar\Delta_{1}(z_{k-1}w_{1})\star\Delta_{1}(z_{l-1}w_{2})\big).
\end{gather*}
using Lemma~\ref{lem:formulas}~(3),~(4) and~(6).
Now apply the induction hypothesis again.
As a~result we f\/ind that $I_{1}(\Delta_{1}(w)\star\Delta_{1}(w'))$ is equal to
\begin{gather*}
x(z_{k-1}w_{1}\shp z_{l}w_{2}+z_{k}w_{1}\shp z_{l-1}w_{2}+\hbar z_{k-1}w_{1}
\shp z_{l-1}w_{2})=z_{k}w_{1}\shp z_{l}w_{2}.
\end{gather*}
This completes the proof.
\end{proof}
\begin{proof}[Proof of Proposition~\ref{prop:star-shuffle}.]
It suf\/f\/ices to prove that $e(w \star w')=e(w)\shp e(w')$ for homogeneous elements $w, w' \in
\widehat{\mathfrak{H}}^{0}$.
If $w=1$ or $w'=1$, then it is trivial.
Now we divide into four cases:
\begin{enumerate}\itemsep=0pt
\item[$(i)$] $w=\xi\rho^{r}w_{1}$ and $w'=\xi\rho^{s}w_{2}$ for $r, s\ge 0$ and $w_{1}, w_{2} \in
\mathfrak{h}^{\ge 1}$,

\item[$(ii)$] $w=z_{k}w_{1}$ and $w'=\xi\rho^{r}w_{2}$ for $k\ge 2$, $w_{1}\in\mathfrak{H}^{1}$, $r\ge 0$ and $w_{2}
\in \mathfrak{h}^{\ge 1}$,

\item[$(iii)$] $w=z_{2}w_{1}$ and $w'=z_{l}w_{2}$ for $l\ge 2$ and $w_{1}, w_{2}\in\mathfrak{H}^{1}$,

\item[$(iv)$] $w=z_{k}w_{1}$ and $w'=z_{l}w_{2}$ for $k, l\ge 3$ and $w_{1}, w_{2}\in\mathfrak{H}^{1}$.
\end{enumerate}

{{\it Case} $(i)$} It holds that
\begin{gather*}
w\star w'=\xi\rho^{r}I_{1}\big(\Delta_{1}(w_{1})\star\Delta_{1}(y\rho^{s}w_{2})\big)
+\xi\rho^{s}I_{1}\big(\Delta_{1}(y\rho^{r}w_{1})\star\Delta_{1}(w_{2})\big)
\\
\phantom{w\star w'=}
{}-\xi\rho^{r+s+1}I_{1}\big(\Delta_{1}(w_{1})\star\Delta_{1}(w_{2})\big).
\end{gather*}
Using Lemma~\ref{lem:main}, we see that
\begin{gather*}
e(w\star w')=\xi\rho^{r}(w_{1}\shp y\rho^{s}w_{2})+\xi\rho^{s}(y\rho^{r}w_{1}\shp w_{2}
)+\xi\rho^{r+s+1}(w_{1}\shp w_{2})=e(w)\shp e(w').
\end{gather*}

{{\it Case} $(ii)$} We proceed by induction on $k$.
Let $k=2$.
Using the result in the Case~$(i)$ and Lemma~\ref{lem:main} we see that
\begin{gather*}
e(w\star w')=eI_{0}e^{-1}\left(e(\xi w_{1})\shp e(\xi\rho^{r}w_{2})\right)+\xi\rho^{r}(z_{2}w_{1}
\shp w_{2}).
\end{gather*}
Because of the equality~\eqref{eq:xishaxi}, it holds that
\begin{gather*}
e(\xi w_{1})\shp e(\xi\rho^{r}w_{2})=\xi\big(w_{1}\shp y\rho^{r}w_{2}+yw_{1}\shp \rho^{r}w_{2}-\rho(w_{1}\shp \rho^{r}w_{2})\big).
\end{gather*}
Using Lemma~\ref{lem:formulas}~(5) we get
\begin{gather*}
eI_{0}e^{-1}\left(e(\xi w_{1})\shp e(\xi\rho^{r}w_{2})\right)
=xy\big(w_{1}\shp y\rho^{r}w_{2}+yw_{1}\shp \rho^{r}w_{2}-\rho(w_{1}\shp \rho^{r}w_{2})\big)
\\
\hphantom{eI_{0}e^{-1}\left(e(\xi w_{1})\shp e(\xi\rho^{r}w_{2})\right)}{}
=x\big(yw_{1}\shp y\rho^{r}w_{2}\big).
\end{gather*}
Hence{\samepage
\begin{gather*}
e(w\star w')=x\big(yw_{1}\shp y\rho^{r}w_{2}\big)+(y-\rho)\rho^{r}(xyw_{1}\shp w_{2}
)=e(w)\shp e(w').
\end{gather*}
This completes the proof for the case $k=2$.}

Suppose that $k\ge 3$.
Using Lemma~\ref{lem:formulas}~(2) we have
\begin{gather*}
e(w\star w')=eI_{0}\big(z_{k-1}w_{1}\star\xi\rho^{r}w_{2}\big)+\xi\rho^{r}I_{1}\big(\Delta_{1}e(z_{k}w_{1}
)\star\Delta_{1}(w_{2})\big).
\end{gather*}
Note that $e(z_{k}w_{1})\in \mathfrak{h}^{\ge 1}$.
From the induction hypothesis and Lemma~\ref{lem:main}, we get
\begin{gather*}
e(w\star w')=eI_{0}e^{-1}(e(z_{k-1}w_{1})\shp \xi\rho^{r}w_{2})+\xi\rho^{r}(e(z_{k}w_{1}
)\shp w_{2}).
\end{gather*}
Because of Lemma~\ref{lem:main}~(1), the second term in the right hand side is equal to
\begin{gather*}
\xi\rho^{r}((x+\hbar)e(z_{k-1}w_{1})\shp w_{2}).
\end{gather*}
Let us calculate the f\/irst term.
Set
\begin{gather*}
\theta_{k}=\sum_{a=2}^{k}\binom{k-2}{a-2}\hbar^{k-a}z_{a-1}\in\sum_{a\ge1}z_{a}\mathfrak{H}^{1}.
\end{gather*}
Then $e(z_{k-1}w_{1})=x\theta_{k-1}w_{1}$.
Hence
\begin{gather*}
e(z_{k-1}w_{1})\shp \xi\rho^{r}w_{2}=x(\theta_{k-1}w_{1}\shp y\rho^{r}w_{2}
)+\xi(x\theta_{k-1}w_{1}\shp \rho^{r}w_{2}).
\end{gather*}
Note that the f\/irst term in the right hand side belongs to $\sum\limits_{a \ge 2}z_{a}\mathfrak{H}^{1}$.
Using Lemma~\ref{lem:main} (5) we f\/ind that
\begin{gather*}
eI_{0}e^{-1}(e(z_{k-1}w_{1})\shp \xi\rho^{r}w_{2})
=(x+\hbar)x(\theta_{k-1}w_{1}\shp y\rho^{r}w_{2})+xy(x\theta_{k-1}w_{1}\shp \rho^{r}w_{2})
\\
\phantom{eI_{0}e^{-1}(e(z_{k-1}w_{1})\shp \xi\rho^{r}w_{2})}
=x\big((x+\hbar)\theta_{k-1}w_{1}\shp y\rho^{r}w_{2}\big).
\end{gather*}
Thus we obtain
\begin{gather*}
e(w\star w')=x\big((x+\hbar)\theta_{k-1}w_{1}\shp y\rho^{r}w_{2}\big)+\xi\rho^{r}
((x+\hbar)x\theta_{k-1}w_{1}\shp w_{2})
\\
\phantom{e(w\star w')}
=x(x+\hbar)\theta_{k-1}w_{1}\shp \xi\rho^{r}w_{2}=e(w)\shp e(w').
\end{gather*}

{{\it Case} $(iii)$} We proceed by induction on $l$.
Let $l=2$.
 From the result in the Case~$(ii)$ we have
\begin{gather*}
e(w\star w')=eI_{0}e^{-1}\left(e(\xi w_{1})\shp e(z_{2}w_{2})+e(z_{2}w_{1})\shp e(\xi w_{2})-\hbar e(\xi w_{1})\shp e(\xi w_{2})\right).
\end{gather*}
It holds that
\begin{gather*}
e(\xi w_{1})\shp e(z_{2}w_{2})+e(z_{2}w_{1})\shp e(\xi w_{2})
-\hbar e(\xi w_{1})\shp e(\xi w_{2})
\\
\qquad
=\xi\left(w_{1}\shp (x-\hbar)yw_{2}+(x-\hbar)yw_{1}\shp w_{2}+\hbar\rho(w_{1}
\shp w_{2})\right)
\\
\qquad\phantom{=}
{}+2xy(w_{1}\shp yw_{2}+yw_{1}\shp w_{2}-\rho(w_{1}\shp w_{2})).
\end{gather*}
Using Lemma~\ref{lem:main}~(5) we get
\begin{gather*}
e(w\star w')=xy\left(w_{1}\shp (x-\hbar)yw_{2}+(x-\hbar)yw_{1}\shp w_{2}
+\hbar\rho(w_{1}\shp w_{2})\right)
\\
\phantom{e(w\star w')=}
{}+2(x+\hbar)xy(w_{1}\shp yw_{2}+yw_{1}\shp w_{2}-\rho(w_{1}\shp w_{2}))
\\
\phantom{e(w\star w')}
=x\left(y(w_{1}\shp xyw_{2}+xyw_{1}\shp w_{2})+2x(yw_{1}\shp yw_{2}
)\right)+\hbar x(yw_{1}\shp yw_{2})
\\
\phantom{e(w\star w')}
=x\left(yw_{1}\shp xyw_{2}+xyw_{1}\shp yw_{2}+\hbar yw_{1}\shp yw_{2}
\right)=e(w)\shp e(w').
\end{gather*}

Next consider the case where $l \ge 3$.
 From the result in the case (ii) and Lemma~\ref{lem:main}~(1), we get
\begin{gather*}
e(w\star w')=eI_{0}e^{-1}(\xi w_{1}\shp xe(z_{l-1}w_{2})+z_{2}w_{1}\shp e(z_{l-1}w_{2})).
\end{gather*}
Note that
\begin{gather*}
\xi w_{1}\shp xe(z_{l-1}w_{2})+z_{2}w_{1}\shp e(z_{l-1}w_{2})
\\
\qquad{}
=\xi\big(w_{1}\shp xe(z_{l-1}w_{2})\big)+x\big(yw_{1}\shp e(z_{l-1}w_{2})\big)+xyw_{1}\shp e(z_{l-1}w_{2}).
\end{gather*}
The second and third terms in the right hand side belong to $\sum\limits_{a \ge 2}z_{a}\mathfrak{H}^{1}$.
Therefore we obtain
\begin{gather*}
e(w\star w')=xy\big(w_{1}\shp xe(z_{l-1}w_{2})\big)
+(x+\hbar)\left\{x\big(yw_{1}\shp e(z_{l-1}w_{2})\big)+xyw_{1}\shp e(z_{l-1}w_{2})\right\}
\\
\phantom{e(w\star w')}{}
=x\big(yw_{1}\shp xe(z_{l-1}w_{2})\big)+\hbar x\big(yw_{1}\shp e(z_{l-1}w_{2})\big)+(x+\hbar)\big(xyw_{1}
\shp e(z_{l-1}w_{2})\big)
\\
\phantom{e(w\star w')}{}
=xyw_{1}\shp xe(z_{l-1}w_{2})+\hbar xyw_{1}\shp e(z_{l-1}w_{2})
\\
\phantom{e(w\star w')}{}
=xyw_{1}\shp (x+\hbar)e(z_{l-1}w_{2})=e(w)\shp e(w')
\end{gather*}
using Lemma~\ref{lem:main}~(1) again.

{{\it Case} $(iv)$} We proceed by induction on $k+l$.
{}From the result in the Case~$(iii)$ and the induction hypothesis we f\/ind that
\begin{gather*}
e(w\star w')\,{=}\,eI_{0}e^{-1}\big(e(z_{k-1}w_{1})\shp e(z_{l}w_{2})
\,{+}\,e(z_{k}w_{1})\shp e(z_{l-1}w_{2})\,{-}\,\hbar e(z_{k-1}w_{1})\shp e(z_{l-1}w_{2})\big).
\end{gather*}
Using Lemma~\ref{lem:main}~(1) we have
\begin{gather*}
e(z_{k-1}w_{1})\shp e(z_{l}w_{2})+e(z_{k}w_{1})\shp e(z_{l-1}w_{2})-\hbar e(z_{k-1}
w_{1})\shp e(z_{l-1}w_{2})
\\
\qquad{}
=e(z_{k-1}w_{1})\shp xe(z_{l-1}w_{2})+xe(z_{k-1}w_{1})\shp e(z_{l-1}w_{2})
+\hbar e(z_{k-1}w_{1})\shp e(z_{l-1}w_{2}).
\end{gather*}
It belongs to $\sum\limits_{a \ge 2}z_{a}\mathfrak{H}^{1}$.
Hence Lemma~\ref{lem:main}~(5) implies that
\begin{gather*}
e(w\star w')=(x+\hbar)\{e(z_{k-1}w_{1})\shp xe(z_{l-1}w_{2})+xe(z_{k-1}w_{1})\shp e(z_{l-1}w_{2})
\\
\phantom{e(w\star w')=}
{}+\hbar e(z_{k-1}w_{1})\shp e(z_{l-1}w_{2})\}
\\
\phantom{e(w\star w')}{}
=(x+\hbar)e(z_{k-1}w_{1})\shp (x+\hbar)e(z_{l-1}w_{2})=e(w)\shp e(w').
\end{gather*}
This completes the proof.
\end{proof}

\section[Linear relations among the modif\/ied $q$MZVs]{Linear relations among the modif\/ied
$\boldsymbol{q}$MZVs}
\label{sec:dimension}

\subsection{Double shuf\/f\/le relation}

We regard $\widehat{\mathfrak{H}}^{0}$ as a~graded $\mathbb{Q}$-module by setting the degree of $x$, $y$, $\rho$
and $\hbar$ to be one, and call the degree the {\it weight} on $\widehat{\mathfrak{H}}^{0}$.
Denote the homogeneous component of weight $d$ by $\widehat{\mathfrak{H}}^{0}_{d}$.
Now we def\/ine the $\mathbb{Q}$-linear map $\bar{Z}_{q}: \widehat{\mathfrak{H}}^{0} \to \mathbb{C}$ by
$\bar{Z}_{q}(w):=(1-q)^{-d}Z_{q}(w)$ for $w \in \widehat{\mathfrak{H}}^{0}_{d}$.
If $\mathbf{k}=(k_{1}, \ldots, k_{r})$ is an admissible index, $\bar{Z}_{q}(z_{k_{1}} \dots z_{k_{r}})$ is
equal to the modif\/ied $q$MZV $\bar{\zeta}_{q}(\mathbf{k})$ def\/ined by~\eqref{eq:def-qMZVbar}.
Set $\mathfrak{H}^{0}_{d}:=\mathfrak{H}^{0}\cap\widehat{\mathfrak{H}}^{0}_{d}$.
Then we have
\begin{gather*}
\bar{Z}_{q}\big(\mathfrak{H}^{0}_{d}\big)=\sum_{\left|{\mathbf{k}}\right|\le d}\mathbb{Q}\bar{\zeta}_{q}
(\mathbf{k}).
\end{gather*}

From the def\/inition of the harmonic product $*$ and the integral shuf\/f\/le product $\mathcyr{sh}$, we
see that $\widehat{\mathfrak{H}}^{0}=\oplus_{d\ge 0}\widehat{\mathfrak{H}}^{0}_{d}$ is a~commutative graded
$\mathbb{Q}$-algebra with respect to either $*$ or $\mathcyr{sh}$.
Now we obtain the following theorem  from Theorems~\ref{thm:harmonic-product} and~\ref{thm:integral-shuffle-product}.
\begin{theorem}
\label{thm:double-shuffle}
Denote by $S_{d}$ $(d\ge 0)$ the $\mathbb{Q}$-subspace of $\widehat{\mathfrak{H}}^{0}_{d}$ spanned by the
elements $w*w'-w\shp w'$ where $w$ and $w'$ are homogeneous elements of $\widehat{\mathfrak{H}}^{0}$
such that the sum of the weights of $w$ and $w'$ is equal to~$d$.
Then $S_{d}\subset \ker{\bar{Z}_{q}}$.
\end{theorem}

Thus we obtain linear relations among the modif\/ied $q$MZVs as the image of $S_{d} \cap \mathfrak{H}^{0}$.
Let us call such relations the {\it double shuffle relations}.

As an example of the double shuf\/f\/le relations we prove a~$q$-analogue of Hof\/fman's identity for
MZVs~\cite{H1}:
\begin{proposition}
\label{prop:Hoffman}
Let $(k_{1}, \ldots, k_{r})$ be an admissible index.
Then we have
\begin{gather*}
\sum_{1\le i\le r}\zeta_{q}(k_{1},\ldots,k_{i}+1,\ldots,k_{r})=\sum_{1\le i\le r\atop{k_{i}\ge2}}\sum_{a=0}
^{k_{i}-2}\zeta_{q}(k_{1},\ldots,k_{i-1},k_{i}-a,a+1,k_{i+1},\ldots,k_{r}).
\end{gather*}
\end{proposition}
\begin{proof}
The proof is similar to that for MZVs given in~\cite{HO}.
{}From the def\/inition of the harmonic product we have
\begin{gather*}
\xi*z_{k_{1}}\cdots z_{k_{r}}=\sum_{i=1}^{r+1}z_{k_{1}}\cdots z_{k_{i-1}}\xi z_{k_{i}}\cdots z_{k_{r}}
+\sum_{i=1}^{r}z_{k_{1}}\cdots z_{k_{i}+1}\cdots z_{k_{r}}.
\end{gather*}
For $\alpha \ge 1$ and $w \in \mathfrak{H}^{1}$, it holds that
\begin{gather*}
y\shp x^{\alpha}w=\sum_{j=0}^{\alpha-1}x^{j}yx^{\alpha-j}w+x^{\alpha}(y\shp w),
\\
y\shp y^{\alpha}w=\sum_{j=1}^{\alpha}y^{j}\xi y^{\alpha-j}w+y^{\alpha}(y\shp w).
\end{gather*}
Using these formulas we obtain
\begin{gather*}
\xi\shp z_{k_{1}}\cdots z_{k_{r}}=\sum_{i=1}^{r+1}z_{k_{1}}\cdots z_{k_{i-1}}\xi z_{k_{i}}
\cdots z_{k_{r}}+\sum_{1\le i\le r\atop{k_{i}\ge2}}\sum_{a=0}^{k_{i}-2}z_{k_{1}}\cdots z_{k_{i-1}}z_{k_{i}
-a}z_{a+1}z_{k_{i+1}}\cdots z_{k_{r}}.
\end{gather*}
Hence we get the desired equality from Theorem~\ref{thm:double-shuffle}.
\end{proof}

\subsection{Resummation duality}

The double shuf\/f\/le relations do not contain all linear relations among the modif\/ied $q$MZVs.
We give another family of linear relations to make up for this lack.
\begin{theorem}
\label{thm:resummation-duality}
For a~positive integer $k$, set
\begin{gather*}
\varphi_{k}:=\sum_{a=2}^{k}(-\hbar)^{k-a}z_{a}+(-\hbar)^{k-1}\xi.
\end{gather*}
Let $r$ be a~positive integer and $\alpha_{i}$, $\beta_{i}$ $(1\le i \le r)$ non-negative integers.
Then we have
\begin{gather}
\bar{Z}_{q}\big(\varphi_{\alpha_{1}+1}\rho^{\beta_{1}}\cdots\varphi_{\alpha_{r}+1}\rho^{\beta_{r}}
\big)=\bar{Z}_{q}(\varphi_{\beta_{r}+1}\rho^{\alpha_{r}}\cdots\varphi_{\beta_{1}+1}\rho^{\alpha_{1}}).
\label{eq:resummation-duality}
\end{gather}
\end{theorem}
\begin{proof}
Note that
\begin{gather*}
I_{\varphi_{k}}(n)=(1-q)^{k}\frac{q^{kn}}{(1-q^{n})^{k}}
\end{gather*}
for $k \ge 1$.
Hence we have
\begin{gather*}
\bar{Z}_{q}\big(\varphi_{\alpha_{1}+1}\rho^{\beta_{1}}\cdots\varphi_{\alpha_{r}+1}\rho^{\beta_{r}}\big)
\\
\qquad
=(1-q)^{\sum\limits_{i=1}^{r}(\alpha_{i}+\beta_{i}+1)}\sum_{n_{1}>\cdots>n_{r}>0}\;\prod_{i=1}^{r}\binom{n_{i}
-n_{i+1}-1}{\beta_{i}}\frac{q^{(\alpha_{i}+1)n_{i}}}{(1-q^{n_{i}})^{\alpha_{i}+1}},
\end{gather*}
where $n_{r+1}=0$.
Expand $1/(1-q^{n_{i}})^{\alpha_{i}+1}$ by using~\eqref{eq:expansion}.
Then we get
\begin{gather*}
(1-q)^{\sum\limits_{i=1}^{r}(\alpha_{i}+\beta_{i}+1)}\sum_{n_{1}>\cdots>n_{r}>0}\sum_{s_{1},\ldots,s_{r}=0}
^{\infty}\;\prod_{i=1}^{r}\binom{n_{i}-n_{i+1}-1}{\beta_{i}}\binom{\alpha_{i}+s_{i}}{\alpha_{i}}q^{(\alpha_{i}
+s_{i}+1)n_{i}}.
\end{gather*}
Now take the sum over $n_{1}, \ldots, n_{r}$ successively using~\eqref{eq:expansion} again.
As a~result we obtain
\begin{gather*}
(1-q)^{\sum\limits_{i=1}^{r}(\alpha_{i}+\beta_{i}+1)}\sum\limits_{s_{1},\ldots,s_{r}=0}^{\infty}\;\prod_{i=1}^{r}
\binom{\alpha_{i}+s_{i}}{\alpha_{i}}\left(\frac{q^{\sum\limits_{j=1}^{i}(\alpha_{j}+s_{j}+1)}}{1-q^{\sum\limits_{j=1}^{i}
(\alpha_{j}+s_{j}+1)}}\right)^{\beta_{i}+1}.
\end{gather*}
Setting $m_{i}=\sum\limits_{j=1}^{i}(\alpha_{i}+s_{i}+1)$ $(1 \le i \le r)$ we see that it is equal to
\begin{gather*}
(1-q)^{\sum\limits_{i=1}^{r}\alpha_{i}}\sum_{m_{r}>\cdots>m_{1}>0}\;\prod_{i=1}^{r}\binom{m_{i}-m_{i-1}-1}{\alpha_{i}
}\frac{q^{(\beta_{i}+1)m_{i}}}{[m_{i}]^{\beta_{i}+1}},
\end{gather*}
where $m_{0}=0$.
This is the right hand side of~\eqref{eq:resummation-duality}.
\end{proof}

Let us call the property~\eqref{eq:resummation-duality} the {\it resummation duality}.
Denote by $R_{d}$ $(d \ge 0)$ the $\mathbb{Q}$-subspace of $\widehat{\mathfrak{H}}^{0}_{d}$ spanned by the
elements
\begin{gather*}
\varphi_{\alpha_{1}+1}\rho^{\beta_{1}}\cdots\varphi_{\alpha_{r}+1}\rho^{\beta_{r}}-\varphi_{\beta_{r}+1}
\rho^{\alpha_{r}}\cdots\varphi_{\beta_{1}+1}\rho^{\alpha_{1}}
\end{gather*}
with $r>0$, $\alpha_{i}, \beta_{i} \ge 0$ $(1\le i \le r)$ and
$\sum\limits_{i=1}^{r}(\alpha_{i}+\beta_{i}+1)=d$.
The resummation duality implies that $R_{d} \subset \ker{\bar{Z}_{q}}$.

Recall that the $\mathbb{Q}$-vector space spanned by the modif\/ied $q$MZVs
\begin{gather*}
Z_{\le d}:=\sum_{\left|\mathbf{k}\right|\le d}\mathbb{Q}\bar{\zeta}_{q}(\mathbf{k})
\end{gather*}
is realized as $\bar{Z}_{q}(\mathfrak{H}^{0}_{d})$ in our framework.
The $\mathbb{Q}$-subspaces $S_{d}$, def\/ined in Theorem~\ref{thm:double-shuffle}, and $R_{d}$ are
contained in $\ker{\bar{Z}_{q}}$.
Therefore the subspace
\begin{gather*}
N_{\le d}:=\mathfrak{H}^{0}\cap(S_{d}+R_{d})
\end{gather*}
describes some linear relations among the modif\/ied $q$MZVs.

By computer experiment we can f\/ind a~lower bound of the dimension of $Z_{\le d}$~\cite{OT}, and calculate
the dimension of $N_{\le d}$.
The result up to weight $7$ is given as follows:
\begin{center}
\begin{tabular}{|c||c|c|c|c|c|c|c|c|c|}
\hline
{}$d$ & 2 & 3 & 4 & 5 & 6 & 7
\\
\hline
{}$\#$ of admissible indices & 1 & 3 & 7 & 15 & 31 & 63
\\
\hline
lower bound of $\dim{Z_{\le d}}$ & 1 & 2 & 4 & 7 & 11 & 18
\\
\hline
{}$\dim{N_{\le d}}$& 0 & 1 & 3 & 8 & 20 & 45
\\
\hline
\end{tabular}
\end{center}
The second line above gives the number of admissible indices whose weight is less than or equal to~$d$.
We see that the sum of the values in the third line and the fourth one is equal to the number of admissible
indices.
Therefore the third line gives the dimension of $Z_{\le d}$ exactly and the space $N_{\le d}$ describes all
linear relations among the modif\/ied $q$MZVs up to weight $7$.

\subsection*{Acknowledgements}

The research of the author is supported by Grant-in-Aid for Young Scientists~(B) No.~23740119.

\pdfbookmark[1]{References}{ref}
\LastPageEnding

\end{document}